\documentclass[12pt]{article}
\usepackage{algorithm2e}

\usepackage[dvips]{graphicx}
\usepackage{textcomp}
\usepackage{amsbsy}
\usepackage{latexsym}
\usepackage[mathscr]{eucal}
\usepackage{amsfonts}
\usepackage{amssymb}

\usepackage[bottom]{footmisc}

\usepackage[usenames,dvipsnames]{xcolor}
\usepackage{tikz,ifthen,fp}
\usepackage{pgfplots}
\usetikzlibrary{backgrounds,arrows}

\usepackage{hyperref}
\usepackage[makeroom]{cancel}

\usepackage{fullpage}

\begin{document}
\bibliographystyle{plain}
\title
{New bounds on the Lebesgue constants of Leja sequences on 
the unit disc and on $\Re$-Leja sequences.}
\author{ 
 Moulay Abdellah CHKIFA 
 \footnote{UPMC Univ Paris 06, UMR 7598, Laboratoire 
 Jacques-Louis Lions, F-75005, Paris, France
\nl
CNRS, UMR 7598, Laboratoire Jacques-Louis Lions, 
F-75005, Paris, France. chkifa@ann.jussieu.fr}
}
\hbadness=10000
\vbadness=10000
\newtheorem{lemma}{Lemma}[section]
\newtheorem{prop}[lemma]{Proposition}
\newtheorem{cor}[lemma]{Corollary}
\newtheorem{theorem}[lemma]{Theorem}
\newtheorem{remark}[lemma]{Remark}
\newtheorem{example}[lemma]{Example}
\newtheorem{definition}[lemma]{Definition}
\newtheorem{proper}[lemma]{Properties}
\newtheorem{assumption}[lemma]{Assumption}
%
\def\RR{\rm \hbox{I\kern-.2em\hbox{R}}}
\def\NN{\rm \hbox{I\kern-.2em\hbox{N}}}
\def\ZZ{\rm {{\rm Z}\kern-.28em{\rm Z}}}
\def\CC{\rm \hbox{C\kern -.5em {\raise .32ex \hbox{$\scriptscriptstyle
|$}}\kern
-.22em{\raise .6ex \hbox{$\scriptscriptstyle |$}}\kern .4em}}
\def\vp{\varphi}
\def\<{\langle}
\def\>{\rangle}
\def\t{\tilde}
\def\i{\infty}
\def\e{\varepsilon}
\def\sm{\setminus}
\def\nl{\newline}
\def\o{\overline}
\def\wt{\widetilde}
\def\wh{\widehat}
\def\cT{{\cal T}}
\def\cA{{\cal A}}
\def\cI{{\cal I}}
\def\cV{{\cal V}}
\def\cB{{\cal B}}
\def\cR{{\cal R}}
\def\cD{{\cal D}}
\def\cP{{\cal P}}
\def\cJ{{\cal J}}
\def\cM{{\cal M}}
\def\cO{{\cal O}}
\def\Chi{\raise .3ex
\hbox{\large $\chi$}} \def\vp{\varphi}
\def\lsima{\hbox{\kern -.6em\raisebox{-1ex}{$~\stackrel{\mboxstyle<}{\sim}~$}}\kern -.4em}
\def\lsim{\hbox{\kern -.2em\raisebox{-1ex}{$~\stackrel{\mboxstyle<}{\sim}~$}}\kern -.2em}
\def\[{\Bigl [}
\def\]{\Bigr ]}
\def\({\Bigl (}
\def\){\Bigr )}
\def\[{\Bigl [}
\def\]{\Bigr ]}
\def\({\Bigl (}
\def\){\Bigr )}
\def\L{\pounds}
\def\pr{{\rm Prob}}
\newcommand{\cs}[1]{{\color{blue}{#1}}}
\def\ds{\displaystyle}
\def\ev#1{\vec{#1}}     
\newcommand{\setu}{{\mathfrak u}}
\newcommand{\lt}{\ell^{2}(\nabla)}
\def\Supp#1{{\rm supp\,}{#1}}
\def\R{\mathbb{R}}
\def\E{\mathbb{E}}
\def\nl{\newline}
\def\T{{\relax\ifmmode I\!\!\hspace{-1pt}T\else$I\!\!\hspace{-1pt}T$\fi}}
\def\N{\mathbb{N}}
\def\L{\mathbb{L}}
\def\P{\mathbb{P}}
\def\V{\mathbb{V}}
\def\D{\mathbb{D}}

\def\Z{\mathbb{Z}}
\def\N{\mathbb{N}}
\def\Zd{\Z^d}
\def\Q{\mathbb{Q}}
\def\C{\mathbb{C}}
\def\Rd{\R^d}
\def\gsim{\mathrel{\raisebox{-4pt}{$\stackrel{\mboxstyle>}{\sim}$}}}
\def\sime{\raisebox{0ex}{$~\stackrel{\mboxstyle\sim}{=}~$}}
\def\lsim{\raisebox{-1ex}{$~\stackrel{\mboxstyle<}{\sim}~$}}
\def\div{{\rm  div }}
\def\M{M}  \def\NN{N}                  
\def\Le{{\ell^1}}            
\def\Lz{{\ell^2}}
\def\Let{{\tilde\ell^1}}     
\def\Lzt{{\tilde\ell^2}}
\def\Ltw{\ell^\tau^w(\nabla)}
\def\t#1{\tilde{#1}}
\def\la{\lambda}
\def\La{\Lambda}
\def\ga{\gamma}
\def\BV{{\rm BV}}
\def\Ga{\eta}
\def\al{\alpha}
\def\cZ{{\cal Z}}
\def\cA{{\cal A}}
\def\cU{{\cal U}}
\def\cV{{\cal V}}
\def\argmin{\mathop{\rm argmin}}
\def\argmax{\mathop{\rm argmax}}
\def\prob{\mathop{\rm prob}}

\def\cO{{\cal O}}
\def\cA{{\cal A}}
\def\cC{{\cal C}}
\def\cF{{\cal F}}
\def\bu{{\bf u}}
\def\bz{{\bf z}}
\def\bZ{{\bf Z}}
\def\bI{{\bf I}}
\def\cE{{\cal E}}
\def\cD{{\cal D}}
\def\cG{{\cal G}}
\def\cI{{\cal I}}
\def\cJ{{\cal J}}
\def\cM{{\cal M}}
\def\cN{{\cal N}}
\def\cT{{\cal T}}
\def\cU{{\cal U}}
\def\cV{{\cal V}}
\def\cW{{\cal W}}
\def\cL{{\cal L}}
\def\cB{{\cal B}}
\def\cG{{\cal G}}
\def\cK{{\cal K}}
\def\cS{{\cal S}}
\def\cP{{\cal P}}
\def\cQ{{\cal Q}}
\def\cR{{\cal R}}
\def\cU{{\cal U}}
\def\bL{{\bf L}}
\def\bl{{\bf l}}
\def\bK{{\bf K}}
\def\bC{{\bf C}}
\def\X{X\in\{L,R\}}
\def\ph{{\varphi}}
\def\H{{\cal H}}
\def\bM{{\bf M}}
\def\bx{{\bf x}}
\def\bj{{\bf j}}
\def\bG{{\bf G}}
\def\bP{{\bf P}}
\def\bW{{\bf W}}
\def\bT{{\bf T}}
\def\bV{{\bf V}}
\def\bv{{\bf v}}
\def\bt{{\bf t}}
\def\bz{{\bf z}}
\def\bw{{\bf w}}
\def \span{{\rm span}}
\def \meas {{\rm meas}}
\def\rhom{{\rho^m}}
\def\lll{\langle}
\def\argmin{\mathop{\rm argmin}}
\def\argmax{\mathop{\rm argmax}}
\def\dJ{\nabla}
\newcommand{\ba}{{\bf a}}
\newcommand{\bb}{{\bf b}}
\newcommand{\bc}{{\bf c}}
\newcommand{\bd}{{\bf d}}
\newcommand{\bs}{{\bf s}}
\newcommand{\bff}{{\bf f}}
\newcommand{\bp}{{\bf p}}
\newcommand{\bg}{{\bf g}}
\newcommand{\by}{{\bf y}}
\newcommand{\br}{{\bf r}}
\newcommand{\be}{\begin{equation}}
\newcommand{\ee}{\end{equation}}
\newcommand{\bea}{$$ \begin{array}{lll}}
\newcommand{\eea}{\end{array} $$}
\def \Vol{\mathop{\rm  Vol}}
\def \mes{\mathop{\rm mes}}
\def \Prob{\mathop{\rm  Prob}}
\def \exp{\mathop{\rm    exp}}
\def \sign{\mathop{\rm   sign}}
\def \sp{\mathop{\rm   span}}
\def \vphi{{\varphi}}
\def \csp{\overline \mathop{\rm   span}}
\newenvironment{disarray}{\everymath{\displaystyle\everymath{}}\array}{\endarray}
\newcommand{\KL}{Karh\'unen-Lo\`eve }
%
\newcommand{\beqn}{\begin{equation}}
\newcommand{\eeqn}{\end{equation}}
\def\beginproof{\noindent{\bf Proof:}~ }
\def\endproof{\hfill\rule{1.5mm}{1.5mm}\\[2mm]}

\newenvironment{Proof}{\noindent{\bf Proof:}\quad}{\endproof}

\renewcommand{\theequation}{\thesection.\arabic{equation}}
\renewcommand{\thefigure}{\thesection.\arabic{figure}}

\makeatletter
\@addtoreset{equation}{section}
\makeatother

\newcommand\abs[1]{\left|#1\right|}
\newcommand\clos{\mathop{\rm clos}\nolimits}
\newcommand\trunc{\mathop{\rm trunc}\nolimits}
\renewcommand\d{d}
\newcommand\dd{d}
\newcommand\diag{\mathop{\rm diag}}
\newcommand\dist{\mathop{\rm dist}}
\newcommand\diam{\mathop{\rm diam}}
\newcommand\cond{\mathop{\rm cond}\nolimits}
\newcommand\eref[1]{{\rm (\ref{#1})}}
\newcommand{\iref}[1]{{\rm (\ref{#1})}}
\newcommand\Hnorm[1]{\norm{#1}_{H^s([0,1])}}
\def\int{\intop\limits}
\renewcommand\labelenumi{(\roman{enumi})}
\newcommand\lnorm[1]{\norm{#1}_{\ell^2(\Z)}}
\newcommand\Lnorm[1]{\norm{#1}_{L^2([0,1])}}
\newcommand\LR{{L^2(\R)}}
\newcommand\LRnorm[1]{\norm{#1}_\LR}
\newcommand\Matrix[2]{\hphantom{#1}_#2#1}
\newcommand\norm[1]{\left\|#1\right\|}
\newcommand\ogauss[1]{\left\lceil#1\right\rceil}
\newcommand{\QED}{\hfill
\raisebox{-2pt}{\rule{5.6pt}{8pt}\rule{4pt}{0pt}}%
  \smallskip\par}
\newcommand\Rscalar[1]{\scalar{#1}_\R}
\newcommand\scalar[1]{\left(#1\right)}
\newcommand\Scalar[1]{\scalar{#1}_{[0,1]}}
\newcommand\Span{\mathop{\rm span}}
\newcommand\supp{\mathop{\rm supp}}
\newcommand\ugauss[1]{\left\lfloor#1\right\rfloor}
\newcommand\with{\, : \,}
\newcommand\Null{{\bf 0}}
\newcommand\bA{{\bf A}}
\newcommand\bB{{\bf B}}
\newcommand\bR{{\bf R}}
\newcommand\bD{{\bf D}}
\newcommand\bE{{\bf E}}
\newcommand\bF{{\bf F}}
\newcommand\bH{{\bf H}}
\newcommand\bU{{\bf U}}
\newcommand\cH{{\cal H}}
\newcommand\sinc{{\rm sinc}}
\def\enorm#1{| \! | \! | #1 | \! | \! |}

\newcommand{\dm}{\frac{d-1}{d}}

\let\bm\bf
\newcommand{\bbeta}{{{\rm \boldmath$\beta$}}}
\newcommand{\bal}{{{\rm \boldmath$\alpha$}}}
\newcommand{\bbi}{{\bm i}}

\def\nnew{\color{black}}
\def\mnew{\color{black}}

\newcommand{\red}[1]{{\color{red}{#1}}}
\newcommand{\blue}[1]{{\color{blue}{#1}}}

\def\substack#1%
    {\begingroup
     \let\\\atop
     #1
     \endgroup}

\newcommand{\dI}{\Delta}
\maketitle
\date{}

%

%
\begin{abstract}
In the papers \cite{Ch1,Ch2}
we have established linear and quadratic bounds, in $k$, on 
the growth of the Lebesgue constants associated with the 
$k$-sections of Leja sequences on the unit disc $\cU$ and 
$\Re$-Leja sequences obtained from the latter by projection 
into $[-1,1]$. In this paper, we improve these bounds and 
derive sub-linear and sub-quadratic bounds. The main
novelty is the introduction of a ``quadratic'' 
Lebesgue function for Leja sequences on $\cU$
which exploits perfectly the binary structure of such sequences 
and can be sharply bounded. This yields new bounds on 
the Lebesgue constants of such sequences, that are almost 
of order $\sqrt k$ when $k$ has a sparse binary expansion. 
It also yields an improvement on the Lebesgue constants
associated with $\Re$-Leja sequences.
\end{abstract}

\section{Introduction}

The growth of the Lebesgue constant of 
Leja sequences on the unit disc and $\Re$-Leja 
sequences was first studied in \cite{CaP1,CaP2}. 
The main motivation was the study of the stability 
of Lagrange interpolation in multi-dimension based 
on intertwining of block unisolvent arrays. Such 
sequences, more particularly $\Re$-Leja sequences, 
were also considered in many other works in the framework 
of structured hierarchical interpolation in high dimension. 
Although not always referred to as such, they are typically 
considered in the framework of sparse grids for interpolation 
and quadrature \cite{HJG,GWG}. Indeed, the sections
of length $2^n+1$ of $\Re$-Leja sequences coincide 
with the Clenshaw-Curtis abscissas $\cos(2^{-n}j\pi),j=0,\dots,2^n$ 
which are de facto used, thanks to the logarithmic growth 
of their Lebesgue constant.

Motivated by the development of cheap and 
stable non-intrusive methods for the treatment 
of parametric PDEs in high dimension, we have 
also used these sequences in \cite{CCS1,CCPP} 
with a highly sparse hierarchical 
polynomial interpolation procedure. The multi-variate 
interpolation process is based on the Smolyak formula and the 
sampling set is progressively enriched in a 
structured infinite grid $\otimes_{j=0}^d Z$
together with the polynomial space 
by only one element at a time.
The Lebesgue constant that quantifies the 
stability of the interpolation process depends 
naturally on the sequence $Z$. We have shown in 
\cite{Ch2} that it has quadratic and cubic 
bounds in the number of points of 
interpolation when $Z$ is a Leja sequence on 
$\cU$ or an $\Re$-Leja sequence, 
thanks to the linear and quadratic bounds on the 
growth of the Lebesugue constant of such sequences,
also established in \cite{Ch1,Ch2}. We refer to the introduction 
and section $2$ in \cite{Ch2} for a concise description 
on the construction of the interpolation process and the 
study of its stability.


The present paper is also concerned with the growth 
of the Lebesgue constant of Leja and $\Re$-Leja 
sequences. We improve the linear and quadratic bounds 
obtained in \cite{Ch2}. In particular, we show that for 
$\Re$-Leja sequnences, the bound is logarithmic 
for many values of $k$ which may be useful for 
proposing cheap and stable interpolation scheme 
in the framework of sparse grids \cite{GWG}.

\subsection{One dimensional hierarchical interpolation}

Let $X$ be a compact domain in $\C$ or $\R$,
typically the complex unit disc $\cU:=\{|z|\leq 1\}$ 
or the unit interval $[-1,1]$, and $Z=(z_j)_{j\geq0}$
a sequence of mutually distinct points in $X$. We 
denote by $I_{Z_k}$ the univariate 
interpolation operator onto the polynomials space 
${\mathbb P}_{k-1}$ associated with the section 
of length $k$, $Z_k:=(z_0,\cdots,z_{k-1})$. The 
interpolation operator is given by Lagrange interpolation 
formula: for $f\in \cC(X)$ and $z\in X$
\be
I_{Z_k} f (z) =\sum_{j=0}^{k-1} f(z_j)  l_{j,k} (z),\quad\quad
l_{j,k}(z):=\prod_{\substack{i=0\\i\neq j}}^{k-1} \frac{z-z_j}{z_i-z_j},
\;\; j=0,\dots,k-1.
\ee
Since the sections $Z_k$ are nested, it is convenient to 
give the operator $I_{Z_k}$ using Newton interpolation 
formula which amounts essentially writing: 
$\Delta_0 f:= I_{Z_1}f \equiv f(z_0)$ and 
\be
\label{newton}
I_{Z_{k+1}}=I_{Z_{k}}+ \Delta_k
= \sum_{l=0}^k \Delta_l, 
\quad\mbox{where}\quad
\Delta_l(Z):=I_{Z_{l+1}}-I_{Z_l},\;\;l\geq1,
\ee
and computing the operators $\Delta_l$ using divided 
differences, see \cite[Chapter 2]{Davis} or equivalently 
the following formula which are differently normalized, see \cite{Ch2,CCS1},
\be
\Delta_l f = \(f(z_l) - I_{Z_l}f(z_l)\)\prod_{j=0}^{l-1} 
\frac {(z-z_j)}{(z_l-z_j)},\quad l\geq1,
\ee

The stability of the operator $I_{Z_k}$ depends on the 
positions of the elements of $Z_k$ on $X$, in particular 
through the Lebesgue constant associated with $Z_k$, 
defined by 
\begin{equation}
\L_{Z_k}:=\max_{f\in C(X)-\{0\}} 
\frac { \|I_{Z_k} f\|_{L^\infty(X)}}{\|f\|_{L^\infty(X)}}
=\max_{z\in X} \lambda_{Z_k}(z),
\label{leblag}
\end{equation}
where $\lambda_{Z_k}$ is the so-called Lebesgue 
function associated with $Z_{k}$ defined by
\be
\label{lagrange}
\lambda_{Z_k}(z):=\sum_{i=0}^{k-1} |l_{i,k}(z)|,
\quad \quad  z\in X.
\ee
We also introduce the notation 
\begin{equation}
\label{normDeltak}
\D_k(Z):=\max_{f\in C(X)-\{0\}} 
\frac { \|\Delta_k f\|_{L^\infty(X)}}{\|f\|_{L^\infty(X)}}.
\end{equation}

In the case of the unit disk or the unit interval, it is known 
that $\L_k$ the Lebesgue constant associated with any 
set of $k$ mutually distinct points can not grow slower 
than $\frac {2\log(k)}\pi$ and it is 
well known that such growth is fulfilled by the set of 
$k$-roots of unity in the case $X=\cU$ and the 
Tchybeshev or Gauss-Lobatto abscissas in the case $X= [-1,1]$, 
see e.g. \cite{Bern}. However such sets of points 
are not the sections of a fixed sequence $Z$.

In \cite{CaP1,CaP2}, the authors considered 
for $Z$ Leja sequences on $\cU$ initiated at the 
boundary $\partial \cU$ and $\Re$-Leja 
sequences obtained by projection onto $[1,1]$  
of the latter when initiated at $1$. They showed 
that the growth of $\L_{Z_k}$ is controlled by ${\cal O}(k\log (k))$ 
and ${\cal O}(k^3\log (k))$ respectively. In our previous works 
\cite{Ch1,Ch2}, we have improved these bounds to $2k$ 
and $8\sqrt 2k^2$ respectively. We have also established 
in \cite{Ch2} the bound ${\mathbb D}_{k}\leq (1+k)^2$ for the 
difference operator, which could not be obtained directly from 
${\mathbb D}_{k}\leq \L_{Z_{k+1}}+\L_{Z_k}$ and which is 
essential to prove that the multivariate interpolation 
operator using $\Re$-Leja sequences has a cubic 
Lebesgue constant, see \cite[formula 25]{Ch2}.

\subsection{Contributions of the paper}
In this paper, we improve the bounds of the
 previous paper 
\cite{CaP1,CaP2,Ch1,Ch2}.
Our techniques of proof share several 
points with those developed in \cite{Ch1,Ch2}, 
yet they are shorter and relies notably on the binary pattern 
of Leja sequences on the unit disk. The novelty 
in the present paper is the introduction of the ``quadractic'' 
Lebesgue constant 
\be
\label{lagrangeQuadratic}
\lambda_{Z_k,2}(z):=
\(\sum_{i=0}^{k-1} |l_{i,k}(z)|^2\)^{\frac12}, \quad z\in X,
\ee
where $l_{i,k}$ are the Lagrange polynomials as defined in \iref{lagrange}.
We study this function and its maximum
\be
\label{lagrangeQuadratic}
\L_{Z_k,2} :=\max_{z\in X} \lambda_{Z_k,2}(z).
\ee
We establish in \S 2 in the case where $Z$ is any Leja sequence on $\cU$ initiated on the boundary
$\partial \cU$ the ``sharp'' inequality
\be
\label{boundLk2Intro}
\lambda_{Z_k}(z_k)
\leq \L_{Z_k,2} \leq 
3\lambda_{Z_k}(z_k),
\quad{and}\quad
\lambda_{Z_k}(z_k):=\sqrt {2^{\sigma_1(k)}-1},
\ee
where $\sigma_1(k)$ denote the number of ones in 
the binary expansion of $k$. Cauchy-schwatrz 
inequality applied to the Lebesgue function 
$\lambda_{Z_k}$ defined in \iref{lagrange} yields 
$\lambda_{Z_k} \leq \sqrt {k}~\lambda_{Z_k,2}$. 
This shows that we also establish
\be
\L_{Z_k} \leq 3 \sqrt k ~\sqrt {2^{\sigma_1(k)}-1},
\ee
for Leja sequences on $\cU$, which improves considerably the 
linear bound $2k$ established in \cite{Ch1} when the binary 
expansion of $k$ is very sparse. For example, for $k=2^n+3$ 
with $n$ large, we get $\L_{Z_k}\leq 3\sqrt{7k} <<k$.
Using the bound \iref{boundLk2Intro}, we establish in \S3 a 
new bound on the growth of Lebesgue constants of $\Re$-Leja 
sequences that implies  
\be
\L_k 
\leq 
6 \sqrt 5
~k~ 2^{\sigma_1(l)}, 
\quad
\mbox{where}
\quad l=k-(2^n+1),
\ee
where $n$ is the integer such that $2^n+1 \leq k <2^{n+1}+1$. Again, we 
remark that the previous bound improves the bound $8\sqrt 2k^2$ established
in \cite{Ch2} when $l$ is small compared to $2^n$ or very sparse in 
the sense of binary expansion. We actually prove a bound
that is logarithmic for many values of $k$ other than the values 
$2^n+1$, see Theorem \ref{theoLebesgueConstant}.

Finally, we provide in \S \ref{sectionDiff} new bounds on the growth 
of $\D_k$ the norm of the difference operators. We provide the bounds 
\be
\D_k \leq 
1+\sqrt{k(2^{\sigma_1(k)}-1)},
\quad \quad 
\D_k \leq 2^{\sigma_1(k)} 2^{n}, \quad k\geq1.
\ee
in the case of Leja sequences on $\cU$ and the case
of $\Re$-Leja sequences respectively where for the latter 
$n$ is defined as above.

\subsection{Notation}

In the remainder of the paper, we work with the following notation. For 
an infinite sequence $Z:=(z_j)_{j\geq0}$ on $X$, we introduce the 
section $Z_{l,m} :=(z_l,\cdots,z_{m-1})$ for any $l \leq m-1$. Given 
two finite sequences $A=(a_0,\dots,a_{k-1})$ and $B=(b_0,\dots,b_{l-1})$, 
we denote by $A\wedge B$ the concatenation of $A$ and $B$, i.e.
$A \wedge B =(a_0,\dots,a_{k-1},b_0,\dots,b_{l-1})$.
For any finite set $S=(s_0,\cdots,s_{k-1})$ of complex numbers
and $\rho \in \C$, we introduce the notation
\begin{equation}
\label{rhoRealConjugate}
\rho S:=(\rho s_0,\cdots, \rho s_{k-1}),\quad 
\Re (S):=(\Re(s_0),\cdots, \Re(s_{k-1})),\quad 
\overline{S}:=(\overline{s_0},\cdots, \overline{s_{k-1}}). 
\end{equation}
Throughout this paper, to any finite set $S$ of numbers,
we associate the polynomial 
\begin{equation}
\label{w}
w_S(z):=\prod_{s\in S}(z-s),\quad 
\mbox{with the convention}
\quad
w_{\emptyset}(z):=1
\end{equation}
Any integer $k\geq1$ can be uniquely expanded according to
\begin{equation}
\label{binaryk}
k = \sum_{j=0}^n a_j 2^j, \quad a_j \in \{0,1\}
\end{equation}
We denote by $\sigma_1(k)$, $\sigma_0(k)$ the number 
of ones and zeros in the binary expansion of $k$ and by 
$p(k)$ 
the largest integer $p$ such that $2^p$ divide $k$. 
For $k=2^n,\dots, 2^{n+1}-1$ with binary expansion 
as above, one has
\be
\sigma_1(k) = \sum_{j=0}^n a_j
\quad\mbox{and}\quad
\sigma_0(k) = \sum_{j=0}^n (1-a_j) = n+1-\sigma_1(k).
\label{defsigma10}
\ee

We should finally note that, unless stated otherwise, we only
work with complex numbers $z$ belonging to the unit circle 
$\partial \cU$. This is because in the complex setting we 
investigate supremums  of sub-harmonic functions, 
$\lambda_{Z_k}$ and $\lambda_{Z_k,2}$, which is always 
attained on the boundary.

\section {Leja sequences on the unit disk}


Leja sequences $E=(e_j)_{j\geq0}$ on $\cal U$ considered in 
\cite{CaP1,CaP2,Ch1,Ch2} 
have all their initial value $e_0\in\partial \cal U$ the unit circle. 
They are defined inductively by picking $e_0 \in \partial{\cal U}$ 
arbitrary and defining $e_k$ for $k\geq1$ by
\begin{equation}
\label{defLejaC}
e_k = {\rm  argmax}_{z\in \cal U} |z-e_{k-1}|\dots |z -e_0|.
\end{equation}
The maximum principle implies that 
$e_j \in \partial \cal U$ for any $j\geq1$. Also, the previous 
${\rm argmax}$ problem might admit many 
solutions and $e_k$ is one of them. We call a $k$-Leja 
section every finite sequence $(e_0,\dots,e_{k-1})$ 
obtained by the same recursive procedure. In particular, 
when $E:=(e_j)_{j\geq1}$ is a Leja sequence 
then the section $E_k = (e_0,\dots,e_{k-1})$ is $k$-Leja section.

In contrast to the interval $[-1,1]$ where Leja sequences 
cannot be computed explicitly,  Leja sequences on 
$\partial \cal U$ are much easier to compute. For instance, 
if $e_0=1$ then we can immediately check 
that $e_1=-1$ and $e_2=\pm i$. Assuming that $e_2=i$ then 
$e_3$ maximizes $|z^2-1||z-i|$, so that $e_3=-i$ because 
$-i$ maximizes jointly $|z^2-1|$ and $|z-i|$. Then $e_4$ 
maximizes $|z^4-1|$, etc. We observe a ``binary patten'' on 
the distribution of the first elements of $E$. 

By radial invariance, an arbitrary 
Leja sequence $E=(e_0,e_1,\dots)$ on $\cU$ with 
$e_0\in\partial\cU$ is merely the product by $e_0$ 
of a Leja sequence with initial value $1$. The latter are 
completely determined according to 
the following theorem, see {\cite{BCCa,CaP1,Ch1}}.
\begin{theorem}
Let $n\geq0$, $2^n < k\leq 2^{n+1}$ and $l=k-2^n$. 
The sequence $E_k=(e_0,\dots,e_{k-1})$, with $e_0=1$,
is a $k$-Leja section if and only if 
$E_{2^n}=(e_0,\dots,e_{2^n-1})$ and 
$U_l= (e_{2^n},\dots,e_{k-1})$
are respectively $2^n$-Leja and $l$-Leja sections 
and $e_{2^n}$ is any $2^n$-root of $-1$.
\label{inversetheostruct}
\end{theorem}

In the light of the previous theorem, a natural construction 
of a Leja sequence $E:=(e_j)_{j\geq0}$ in $\cU$ follows
by the recursion
\begin{equation}
E_1 := (e_0=1)
\quad  {\rm and}\quad 
E_{2^{n+1}} :=E_{2^n} \wedge e^{\frac {i\pi}{2^n}} E_{2^n},\quad n\geq0.
\end{equation}
This recursive construction of the sequence $E$ yields an interesting distribution 
of its elements. Indeed, by an immediate induction, see \cite{BCCa}, it can be 
shown that the elements $e_k$ are given by
\begin{equation}
\label{SimpleLejaSequence}
e_k 
=
 \exp
\Big(
i\pi \sum_{j=0}^n a_j 2^{-j}
\Big) 
\quad {\rm for}\quad
 k = \sum_{j=0}^n a_j 2^j, \quad a_j \in \{0,1\}.
\end{equation}
The construction yields then a low-discrepancy sequence on 
$\partial \cal U$ based on the bit-reversal Van der Corput 
enumeration.

As already mentioned above, Theorem \ref{inversetheostruct} 
characterizes completely Leja sequences on the unit circle. It 
has also many implications that turn out to be very useful in 
the analysis of the growth of Lebesgue constants. 
\begin{theorem}
Let $E:=(e_j)_{j\geq0}$ be a Leja sequence on $\cal U$ 
initiated at $e_0\in \partial \cU$. We have:
\begin{itemize}
\item For any $n\geq0$, $E_{2^n}= e_0 {\cal U}_{2^n}$ in the set
sense where ${\cal U}_{2^n}$ is the set of $2^n$-root of unity.
\item For any $k\geq1$, $|w_{E_k}(e_k)| =\sup_{z\in\partial{\cal U}} |w_{E_k}(z)| = 2^{\sigma_1(k)}$.
\item For any $n\geq0$, $E_{2^n,2^{n+1}} := (e_{2^n},\cdots,e_{2^{n+1}-1})$ is a $2^n$-Leja section.
\item The sequence $E^2:=(e_{2j}^2)_{j\geq0}$ is 
a Leja sequence on $\partial \cU$.
\end{itemize}
\label{TheoImplications}
\end{theorem}

Such properties can be easily checked for the simple 
sequence defined in \iref{SimpleLejaSequence}
and are given in \cite{CaP1,Ch1} for more general
Leja sequences.

\subsection{Analysis of the quadratic Lebesgue function}

It is proved in \cite{Ch1} that given two $k$-Leja sections $E_k$ and $F_k$, 
one has $F_k=\rho E_k$ in the set sense for some $\rho\in \partial \cU$. This 
means that the sequence $F_k$ can be obtained from $E_k$ by a permutation 
and the product by $\rho$. By inspection of the quadratic Lebesgue function
\iref{lagrangeQuadratic}, we have then that 
\be
\lambda_{F_k,2}(z) = \lambda_{E_k,2}(z/\rho),\;\;\;z\in \cU 
\quad
\Longrightarrow
\quad
\L_{F_k,2} = \L_{E_k,2}.
\ee
In order to compute the growth of $\L_{E_k,2}$ for arbitrary Leja 
sequences $E$, it suffices then to consider $E$ to be the 
simple sequence given by \iref{SimpleLejaSequence}. Unless 
stated otherwise, for the rest of this section, $E$ is exclusively 
used for this notation. Let us note that 
\be
\label{memoryShort}
E^2 :=(e_{2j}^2)_{j\geq0}=E.
\ee

In order to study the functions  $\lambda_{E_k,2}$, we adopt the 
methodology that we introduced in \cite{Ch1}. Namely, we study 
the implication of $E$ being a Leja sequence in general, on the 
growth of $\lambda_{E_k,2}$, then we use the implication of the 
particular binary distribution of $E$ to derive such growth.

\begin{lemma}
Let $Z$ be a Leja sequence on a real or complex compact $X$. For 
any $k\geq 1$ and any $z\in X$, it holds
\be
\label{growthLeja}
\lambda_{Z_{k+1},2}(z) 
\leq
\lambda_{Z_{k},2}(z)
+
\lambda_{Z_{k},2}(z_k)+1.
\label{kkp1}
\ee
\label{propertyLeja2}
\end{lemma}
{\bf Proof:}
We fix $k\geq 1$ and denote by $l_0,\dots,l_{k-1}$ the Lagrange polynomials 
associated with the section $Z_{k}$ and by $L_0,\dots,L_{k}$ the Lagrange 
polynomials associated with the section $Z_{k+1}$. By Lagrange interpolation 
formula, for $j=0,\dots,k-1$ 
$$
l_j = \sum_{i=0}^k l_{j}(z_i) L_i = L_j+ l_j(z_k) L_k
\quad\Rightarrow\quad
L_j= l_j- l_j(z_k) L_k.
$$
We have then for any $z\in X$
$$
\(\sum_{j=0}^{k-1}|L_j (z)|^2\)^{1/2}
\leq\(\sum_{j=0}^{k-1}|l_j(z)|^2\)^{1/2}+ 
|L_k(z)|\(\sum_{j=0}^{k-1} |l_j(z_k)|^2\)^{1/2}.
$$
where we have merely applied triangular 
inequality with the euclidean norm in $\C^k$.
This also writes
$$
\(|\lambda_{Z_{k+1},2}(z)|^2 
- |L_k(z)|^2\)^{\frac 12}
\leq
\lambda_{Z_k,2}(z) 
+
|L_k(z)|
\lambda_{Z_k,2}(z_k). 
$$
We conclude the proof using $a \leq \sqrt{a^2-b^2}+b$ for $a\geq b\geq0$,
and the inequality
$$
|L_k(z)| = \frac{|w_{Z_k}(z)|}{|w_{Z_k}(z_k)|} \leq 1,
$$
which follows from the Leja definition \iref{defLejaC}.\hfill $\blacksquare$
\nl

The previous result shows that given $Z$ a Leja sequence 
over $X$, the growth of $\L_{Z_k,2}$ is monitored by the 
growth of $\lambda_{Z_k,2}(z_k)$. In particular, it is easily 
checked using induction on $k$ that
\be
\lambda_{Z_k,2}(z_k) = \cO(\log(k))\quad \Longrightarrow \quad \L_{Z_k,2} =  \cO(k\log(k)),
\ee
and
\be
\lambda_{Z_k,2}(z_k) = \cO(k^\theta) \quad \Longrightarrow \quad \L_{Z_k,2} =  \cO(k^{\theta+1}).
\ee

In the following, we show basically that the previous 
implication holds with $\theta=1/2$ for Leja sequences on $\cU$.
However, we use the particular structure of such 
sequences in order to show that the exponent 
$\theta=1/2$ is not deteriorated and that it is also valid 
for $\L_{E_k,2}$. We recall that we work with the simple 
sequence $E$ given in \iref{SimpleLejaSequence} for 
which $E^2=E$. The binary patten of the distribution of E 
on the unit disc yields the following result.

\begin{lemma}
Let $E$ be as in \iref{SimpleLejaSequence}. For any $N\geq1$, one has
 \be
 \lambda_{E_{2N},2}(z) = \lambda_{E_N,2} (z^2),
 \quad \quad z\in \partial \cU.
 \label{k2nk}
\ee
\label{TheoBoundRecurs}
\end{lemma}
{\bf Proof:}
Let $l_0,\dots,l_{2N-1}$ be the Lagrange polynomials associated 
with $E_{2N}$ and $L_0,\dots,L_{N-1}$ be the Lagrange polynomials 
associated with $E_N$. Since $e_{2j+1}= -e_{2j}$ for any $j\geq0$, 
then in view of \iref{memoryShort}
$$
w_{E_{2N}}(z) = w_{E_N^2}(z^2) = w_{E_N}(z^2).
$$
Deriving with respect to $z$ and using 
$(e_{2j+1})^2 =(e_{2j})^2= e_j$ for any $j\geq0$, we deduce that 
\be
\label{Deriv}
|w_{E_{2N}}'(e_{2j+1})|
=
|w_{E_{2N}}'(e_{2j})|
=
2 |w_{E_N}'(e_{2j}^2)|
=
2 |w_{E_N}'(e_j)|,\quad j\geq0.
\ee
We have for any $j=0,\dots,N-1$
$$
|l_{2j}(z)|
=\frac {|w_{E_{2N}}(z)|}{|w_{E_{2N}}'(e_{2j})| |z-{e_{2j}}|},
\quad
\quad 
\quad
|l_{2j+1}(z)|
=\frac {|w_{E_{2N}}(z)|}{|w_{E_{2N}}'(e_{2j+1})| |z-{e_{2j+1}}|}.
$$
Therefore in view of the previous equalities 
\be
|l_{2j}(z)|^2
+
|l_{2j+1}(z)|^2
=
\frac {|w_{E_N}(z^2)|^2}{4 |w_{E_N}'(e_j)|^2 }
\Big[
\frac 1{|z-e_{2j}|^2}
+
\frac {1}{|z+{e_{2j}}|^2}
\Big]
=
\frac {|w_{E_N}(z^2)|^2}{ |w_{E_N}'(e_j)|^2 |z^2-e_j|^2}
=
|L_j(z^2)|^2,
\label{sumTwoPolynomials}
\ee
where we have used $|a-b|^2 +|a+b|^2 = 4$ for 
$a,b\in\partial \cU$ and $e_{2j}^2=e_j$. Summing the 
previous identities for the indices $j=0,\dots,N-1$, 
we get the result. \hfill $\blacksquare$
\nl

We note that the previous result combined with 
$E_{2^n}=\cU_{2^n}$ in the set sense implies that 
\be
\label{lambda2E2n}
\sum_{j=0}^{2^n-1} 
\Big|
\frac{z^{2^n}-1}{2^n(z-e_j)}
\Big|^2
=\lambda_{E_{2^n},2}(z)
=\lambda_{E_1,2}(z^{2^n})=1,
\ee
for any $z\in \partial \cU$. We now turn to the 
growth of $\lambda_{E_k,2}(e_k)$, which 
as mentioned earlier monitor the growth of $\L_{E_k,2}$.

\begin{lemma}
For the Leja sequence $E$  defined in 
\iref{SimpleLejaSequence}, we have for any $k\geq1$, 
\be
\label{lambdaEkk}
\lambda_{E_k,2}(e_k) = \sqrt{ 2^{\sigma_1(k)}-1}
\ee
\label{lemmaQuadraticEkek}
\end{lemma}
{\bf Proof:} First, by Lemma \ref{TheoBoundRecurs} and 
$e_{2N}^2=e_N$, one has 
\be
\label{k2nkAtek}
|\lambda_{E_{2N},2}(e_{2N})|^2 = 
|\lambda_{E_N,2} (e_N)|^2,\quad N\geq1.
\ee
Let now $k$ be an odd number and we write $k=2N+1$ with $N\geq1$. Let 
$l_0,\dots,l_{2N}$ be the Lagrange polynomials associated with 
$E_k$ and $L_0,\dots,L_{N-1}$ be the Lagrange polynomials associated 
with $E_N$. For any $m=0,\dots,2N$, one has 
$$
l_m(e_k) = 
\frac {w_{E_k}(e_k)}{(e_k-e_m) w_{E_k}'(e_m)}
=
\frac {w'_{E_{k+1}}(e_k)}{w_{E_{k+1}}'(e_m)}
\quad\Rightarrow\quad
|l_m(e_k)| = 
\frac {|w'_{E_{N+1}}(e_k^2)|}{|w_{E_{N+1}}'(e_m^2)|},
$$
where we have used $k+1=2(N+1)$ and \iref{Deriv}.
Using $e_{k}^2=e_N$ and $(e_{2j+1})^2 =(e_{2j})^2= e_j$ for any $j$, 
we get for 
$m=2j$ or $m=2j+1$ with $j=0,\dots,N-1$ 
$$
|l_m(e_k)| 
= \frac {|w'_{E_{N+1}}(e_N)|}{|w_{E_{N+1}}'(e_j)|}
=|L_j(e_N)|
\quad\mbox{and also}\quad 
|l_{2N}(e_k)| = \frac {|w'_{E_{N+1}}(e_N)|}{|w_{E_{N+1}}'(e_N)|}=1. 
$$
Summing the numbers $|l_m(e_k)|^2$ over $m=0,\dots,2N$, we infer 
\be
\label{recursionEkek}
|\lambda_{E_{2N+1},2}(e_{2N+1})|^2 = 2 |\lambda_{E_N,2}(e_N)|^2+1.
\ee
In view of the above and $\lambda_{E_1,2}(e_1)$=1, the sequence 
$\alpha:=(\alpha_k := |\lambda_{E_k,2}(e_k)|^2)_{k\geq1}$
satisfies:
$$
\alpha_1=1\quad \mbox{and}\quad
\alpha_{2N}=\alpha_N,\quad \alpha_{2N+1}=2\alpha_N+1, \quad N\geq1.
$$
We have 
$\sigma_1(1)=1$ and $\sigma_1(2N)=\sigma_1(N)$, 
$\sigma_1(2N+1)=\sigma_1(N)+1$ for any $N\geq1$.
It is then easily checked that $(2^{\sigma_1(k)}-1)_{k\geq1}$
satisfies the same recursion as $\alpha$. This shows that 
$\alpha_k=2^{\sigma_1(k)}-1$ for any $k\geq1$ and finishes the proof.
\hfill $\blacksquare$
\nl

We are now able to conclude the main result of this section, which 
states basically that for the sequence $E$ or more generally any 
Leja sequence on $\cU$ initiated at the boundary $\partial\cU$, the 
value of $\L_{E_k,2}=\max_{z\in\cU}\lambda_{E_k,2}(z)$ is 
almost equal to $\lambda_{E_k,2}(e_k)$.
\begin{theorem}
For the Leja sequence $E$  defined in 
\iref{SimpleLejaSequence}, we have for 
any $k\geq1$ 
\be
\label{boundLinearLambdak}
1\leq \frac {\L_{E_k,2}}{\lambda_{E_k,2}(e_k)} 
=\frac {\L_{E_k,2}}{\sqrt{ 2^{\sigma_1(k)}-1}} \leq 3
\ee
\label{theoremLEk2}
\end{theorem}
{\bf Proof:} 
The first part of the inequality is immediate from the definition 
of $\L_{E_k,2}$. Also in view Lemma \ref{TheoBoundRecurs}
and formula \iref{k2nkAtek}, we only need to show 
\iref{boundLinearLambdak} when $k$ is an odd number. Let 
$k=2N+1$ with $N\geq1$. Using 
Lemma \ref{propertyLeja2}, Lemma \ref{TheoBoundRecurs}
and formula \iref{k2nkAtek}, we have 
$$
\lambda_{E_k,2}(z)
\leq \lambda_{E_{2N},2}(z)
+\lambda_{E_{2N},2}(e_{2N})+1
=\lambda_{E_N,2}(z^2)+
\lambda_{E_N,2}(e_N)+1.
$$
If we assumes that 
$\lambda_{E_N,2}(z^2)\leq 3 \lambda_{E_N,2}(e_N)$, we get
$$
\lambda_{E_k,2}(z)
\leq 4\lambda_{E_N,2}(e_N)+1
\leq
3 \sqrt{2|\lambda_{E_N,2}(e_N)|^2+1},
$$
where we have used the elementary inequality $4t+1\leq 3\sqrt{2t^2+1}$
for any $t\geq0$. In view of \iref{recursionEkek}, one then gets 
$\lambda_{E_k,2}(z)\leq 3\lambda_{E_k,2}(e_k)$. The 
verification $\L_{E_1,2}= \lambda_{E_1,2}(e_1)=1$ shows that 
the result follows using an induction on $k\geq1$.\hfill $\blacksquare$

\subsection{Implications on the Lebesgue constant}

The methodology we have provided so far for bounding $\L_{E_k,2}$
is not new, we have developed it in \cite{Ch1} in order to give linear
estimate for $\L_{E_k}$, namely $\L_{E_k}\leq 2k$. Theorem
\iref{theoremLEk2} has also implications on the growth of the Lebesgue 
constant $\L_{E_k}$. Indeed, Cauchy Schwartz inequality applied to the 
Lebesgue function $\lambda_{E_k}$ implies 
$\lambda_{E_k}\leq~ \sqrt k~ \lambda_{E_k,2}$, so that 
\be
\label{boundLkLatest}
\L_{E_k} 
\leq \sqrt k~ \L_{E_k,2}
\leq 3\sqrt {k(2^{\sigma_1(k)}-1)}. 
\ee
The Cauchy Schwartz formula
$\lambda_{E_k}\leq \sqrt {k} ~\lambda_{E_{k,2}}$
is possibly not very pessimistic. It has been recently 
proved that the Lagrange polynomials are uniformly 
bounded, see \cite{Irigoyen}
We shall observe in particular, see Figure, 
that the binary pattern observed for the exact value 
of $\L_{E_k}$ is captured by the previous bound. Moreover,
we are able to provide a lower bound for $\L_{E_k}$,
that is comparable to the previous upper bound for values of 
$k$ with full binary expansion.

\begin{prop}
For the Leja sequence $E$  defined in 
\iref{SimpleLejaSequence}, we have for 
any $k\geq1$ 
\label{TheoremSharpnessLebesgue}
\be
2^{\sigma_1(k)}-1
\leq
\lambda_{E_k}(e_k)
\leq
\L_{E_k}.
\ee
\end{prop}
{\bf Proof:}
We let $N\geq1$ and we use the notation of the proof of Lemma 
\ref{TheoBoundRecurs}. As for formula 
\iref{sumTwoPolynomials} and since 
$|a-b| +|a+b| \geq2$ for 
any $a,b\in\partial \cU$, one has
$$
|l_{2j}(z)|
+
|l_{2j+1}(z)|
=
\frac {|w_{E_N}(z^2)|}{2 |w_{E_N}'(e_j)| }
\frac {|z-e_{2j}|+|z+e_{2j}|}{|z-e_j|}
\geq 
|L_j(z^2)|.
$$
This implies 
$\lambda_{E_{2N}}(z)\geq \lambda_{E_N}(z^2)$ and more particularly 
$\lambda_{E_{2N}}(e_{2N})\geq \lambda_{E_N}(e_N)$.
As in the proof of Lemma \ref{lemmaQuadraticEkek}, we have 
also $\lambda_{E_{2N+1}}(e_{2N+1}) =  2\lambda_{E_{N}}(e_N)+1$.
The sequence 
$(b_k := \lambda_{E_k}(e_k))_{k\geq1}$
satisfies: 
$$
b_1=1\quad \mbox{and}\quad
b_{2N}\geq b_N,\quad b_{2N+1}=2b_N+1, \quad N\geq1.
$$
The sequence $b$ then satisfies $b_k \geq 2^{\sigma_1(k)}-1$ 
for any $k\geq1$.
\hfill $\blacksquare$
\nl 

The previous theorem combined with Theorem \ref{theoremLEk2}
and \iref{boundLkLatest} implies 
\be
\frac{\sqrt{2^{\sigma_1(k)}-1}} 3 
\L_{E_k,2}
\leq
\L_{E_k} 
\leq \sqrt k~ \L_{E_k,2}.
\ee
Cauchy Schwartz inequality is then satisfactory when 
$k\simeq 2^{\sigma_1(k)}$, that is when $k$ has a full 
binary expansion.
\begin{remark}
For integers $k=2^n,\dots,2^{n+1}-1$, if $k=2^{n+1}-1$ in 
which case $\sigma_1(k)=n+1$ is the largest possible, the 
bound \iref{boundLkLatest} merely implies $\L_{E_k}\leq 3k$ which is 
worse than the bound $2k$ established \cite{Ch1} and the 
exact value $\L_{E_k}=k$ of this case, see \cite{CaP1}. 
However, since $\sigma_1(k)=n+1-\sigma_0(k)$ for any  
$k =2^n,\dots 2^{n+1}-1$, then by \iref{boundLkLatest}
\be
\L_{E_k} 
\leq \sqrt{\frac{18}{2^{\sigma_0(k)}}}~\sqrt {2^nk} 
\leq \sqrt{\frac{18}{2^{\sigma_0(k)}}}~k.
\ee
This shows in particular that $\L_{E_k} \leq k$ 
whenever $\sigma_0(k)\geq5$. This last result 
answers partly the conjecture raised in \cite{CaP1}
and which states that $\L_{E_k} \leq k$ for any $k\geq1$.
\end{remark}

For the purpose of the next section, we improve
the bound \iref{boundLkLatest} in the case where 
$k$ is an even number. We recall that we have shown 
in \cite[Theorem 2.8]{Ch1}
\be
\label{binarymemory}
\L_{E_{2^p l}}\leq \L_{2^p}\L_{E_l},
\quad
\quad 
p\geq0,\;\; l \geq1,
\ee
where $\L_{2^p}$ is the Lebesgue 
constant associated with the set of $2^p$-roots of unity. 
The value $\L_{2^p}$ can be computed easily for small 
values of $p$ and it grows logarithmically in $2^p$, 
see e.g. \cite[formula 2.25]{Ch1},
\be
\label{boundL2p}
\L_1=1,\quad\L_2=\sqrt 2,\quad \mbox{and}\quad 
\L_{2^p}\leq  \frac 2\pi \(\log(2^p)+ 9/4\),\quad p\geq2.
\ee
Since $\sigma_1(k)=\sigma_1(k/2^{p(k)})$, we have 
then in view of \iref{boundLkLatest}
and \iref{binarymemory} the following 
theorem

\begin{theorem}
Let $E$ be the Leja sequence defined in 
\iref{SimpleLejaSequence} or any Leja 
sequence on $\cU$ initiated at $\partial \cU$. 
We have
\be
\label{bestBound}
\L_{E_k} 
\leq 3\sqrt {\frac {k}{2^{p(k)}}(2^{\sigma_1(k)}-1)}
~~\L_{2^{p(k)}},\quad k\geq1.
\ee
\label{theoremBoundLeja}
\end{theorem}

We should mention that our primary interest in studying 
$\lambda_{E_k,2}$ was the improvement of the 
results of \cite{Ch2} concerned with the Lebesgue 
constants of $\Re$-Leja sequences. This will be made 
clear in the proof of Theroem \ref{theoLebesgueConstant}. For the 
sake of the same theorem, we need also to provide a 
growth property of Leja sequences on the unit disc.

We let $E=(e_j)_{j\geq0}$ be the simple Leja 
sequence defined by \iref{SimpleLejaSequence}.
For $m\geq0$ and $1\leq l\leq 2^{m-1}$, we 
introduce the notation $K=2^m+l$ and 
$F_{m,l} = E_{2^m,K}$ 
and define the quantity 
\be
\label{defBetaln}
\gamma_{m,l} =\frac 1 {4^m} 
\sum_{j=0}^{K-1} 
\frac {4}{|w_{F_{m,l}}(\overline{e_j})|^2}. 
\ee
The quantity $\gamma_{m,l} $ is well defined. Indeed, by
the particular structure of the sequence $E$, we have 
$E_{2^m+2^{m-1}}=E_{2^m}\wedge e^{\frac {i\pi}{2^m}} E_{2^{m-1}}$, 
so that $E_{2^m+2^{m-1}}=\cU_{2^m}\wedge e^{\frac {i\pi}{2^m}}~\cU_{2^{m-1}}$ in the 
set sense. We have then for $j=0,\dots,2^m+l-1$,
$\o{e_j}$ is in $\cU_{2^m}\wedge e^{\frac {-i\pi}{2^m}}~\cU_{2^{m-1}}$
which does not intersect with $F_{m,l}\subset e^{\frac {i\pi}{2^m}}~\cU_{2^{m-1}}$. We have the
following growth for $\gamma_{m,l}$. 
\begin{lemma}
For any $m\geq1$ and any $1\leq l\leq 2^{m-1}$, we have 
\be
\gamma_{m,l} 
\leq 
\frac {5}{2^{\sigma_1(l)+p(l)+1}}
\ee
\label{lemmabetanl}
\end{lemma}
{\bf Proof:} Since $(e_0,e_1,e_2)=(1,-1, i)$, it can be 
checked that $\gamma_{1,1}=5/4$. We then fix $m\geq2$. 
We define $\rho = e_{2^m}= e^{i\pi/2^m}$, so that 
$F_{m,1} = \{\rho\}$. We have 
$$
\gamma_{m,1}= 
\sum_{j=0}^{2^m} 
\frac {4}{(2^m|e_j - \overline{\rho}|)^2}
=
|\lambda_{E_{2^m},2} (\overline{\rho})|^2
+
\frac {4}{(2^m|\rho - \overline{\rho}|)^2} 
=
1
+\frac 1{|2^m \sin(\pi/2^m)|^2} 
$$
where we have used \iref{lambda2E2n} 
and used that $\o\rho$ is a $2^m$-root of $-1$.
Since
$2^m\sin(\pi/{2^m})\geq 2$ then $\gamma_{m,1} \leq 5/4$. 
For the other values of $l=2,\dots,2^{m-1}$, we have
\begin{itemize}
\item If $l=2N$, we have for any $j\geq0$ that 
$w_{F_{m,l}}(\overline{e_{2j+1}})=w_{F_{m,l}}(\overline{e_{2j}})
=w_{E_{2^{m-1},2^{m-1}+N}}(\overline{e_j})$.
Pairing the indices in \iref{defBetaln} as $2j$ and $2j+1$ with 
$j=0,\dots,2^{n-1}+N-1$, we deduce 
$$
\gamma_{m,l} =\frac{ \gamma_{m-1,N}}{2}.
$$

\item If $l=2N+1$ with $N\geq1$, we may write 
$$
\gamma_{m,l} 
=
\frac 1 {4^{m-1}}
\sum_{j=0}^{K-1} 
\frac {|e_K-\overline{e_j}|^2}{|w_{F_{m,l+1}}(\overline{e_j})|^2}
\leq
\frac 1 {4^{m-1}}
\sum_{j=0}^{K} 
\frac {|e_K-\overline{e_j}|^2}{|w_{F_{m,l+1}}(\overline{e_j})|^2}
=
\gamma_{m-1,N+1},
$$
where we have again paired the indices by $2j$ and $2j+1$ 
for $j=0,\dots,2^n+(N+1)-1$ and used $e_{2j+1}=-e_{2j}$ and 
the identity $|a+b|^2+|a-b|^2=4$ for any $a,b\in\partial \cU$. 
\end{itemize}
Therefore 
$$
\gamma_{m,l} \leq \frac{~5~}4 a_{m,l},
\quad\quad 1\leq m,\;\; 1\leq l\leq2^{m-1},
$$
where $(a_{m,l})_{\substack{1\leq m\\1\leq l\leq2^{m-1}}}$ is the 
sequence that saturates the previous inequalities and hence 
is defined by the following recursion: 
$$
a_{m,1}=1,\;\;m\geq1
\quad\mbox{and}\quad
\left\{
\begin{array}{l}
a_{m,2N}\;\;\;=a_{m-1,N}{\big/} 2\quad\quad n\geq1, N=1,\dots,2^{m-2},\\
a_{m,2N+1}=a_{m-1,N+1}\quad\quad n\geq1, N=1,\dots,2^{m-2}-1.
\end{array}
\right.
$$
The sequence $(a_{m,l})$ has no dependance on $m$ and 
it is equal, in the sense $a_{m,l}= a_l$, to the sequence 
$(a_l)_{l\geq1}$ which 
satisfies the recursion: 
$a_1=1,\;a_{2N}=a_N/2,\; a_{2N+1}=a_{N+1}$.
Since  $\sigma_1(1) + p(1)=1$, 
$\sigma_1(2N) + p(2N)=\sigma_1(N)+p(N)+1$
and
$$
\sigma_1(2N+1) + p(2N+1)
=\sigma_1(2N+1) 
=\sigma_1(N) +1
=\sigma_1(N+1) +p(N+1),
$$
then an immediate induction 
shows that $a_l=2^{1-\sigma_1(l)-p(l)}$,
which finishes the proof. \hfill $\blacksquare$

\section{$\Re$-Leja sequences on $[-1,1]$}

$\Re$-Leja sequences were introduced and studied in 
\cite{CaP2}. Such sequences are simply defined as the 
projection, element-wise but without repetition, into [-1,1] of Leja sequences 
on $\cU$ initiated at $1$. More precisely, given $E=(e_j)_{j\geq0}$ a 
Leja sequence on $\cU$ initiated at $1$, the $\Re$-Leja 
sequence $R=(r_j)_{j\geq0}$ associated with $E$ is obtained 
progressively by: $r_0=\Re(e_0)=1$, $J(0)=0$ and 
\be
\label{algoRLeja}
r_k= \Re(e_{J(k)}) 
\quad\mbox{where}\quad 
J(k)=\min\{j>J(k-1):\Re(e_j)\not\in R_k\},\quad k\geq1.
\ee
This means one projects $e_j$ if and only if $e_j\neq\o{e_i}$
for all $i<j$. The projection rule that prevents the repetition is 
provided in \cite[Theorem 2.4]{CaP2}. One has 
\be
R = \Re (\Xi),
\quad {\rm with } \quad 
\Xi:=(1,-1) \wedge \bigwedge_{j=1}^\infty E_{2^j,2^j+2^{j-1}}.
\label{formRLeja}
\ee
Using a simple cardinality argument, see  \cite[Theorem 2.4]{CaP2} or 
\cite[Formula 40]{Ch2}, this implies 
that the function $J$ used in \iref{algoRLeja} is given by: 
$J(0)=0,~J(1)=1$ and 
\be
\label{formulaJ}
J(k) = 2^n+k-1, \quad\quad n\geq0, \;\;\; 2^n+1 \leq k <2^{n+1}+1.
\ee
In view of \iref{formRLeja} and the properties of Leja sequences 
on $\cU$, any $\Re$-Leja sequence $R$ satisfies  
$r_0=1,~r_1=-1,~r_2=0$ and $r_{2j-1} = -r_{2j}$ for any $j\geq2$. 
An accessible example of an $\Re$-Leja sequence is the one
associated with the simple Leja sequence given by the bit-reversal 
enumeration \iref{SimpleLejaSequence}. We have shown  
in \cite{Ch1} that $R=(\cos(\phi_j))_{j\geq0}$ where the sequence of 
angles $(\phi_k)_{k\geq0}$ is defined recursively by $\phi_0=0$, $\phi_1=\pi$, 
$\phi_2 = \pi/2$ and
\be
\phi_{2j-1} = \frac{\phi_j}{2},\quad\quad
\phi_{2j} = \phi_{2j-1} +\pi,\;\;\;j\geq2.
\ee
This recursion provides a simple process to 
compute an $\Re$-Leja sequence.
We can also construct a Leja sequence 
by simply using the recursion 
$r_0=1$, $r_1=-1$, $r_2=0$ and 
\be
r_{2j-1} = \sqrt{\frac{r_j+1}{2}},\quad\quad
r_{2j} = - r_{2j-1},\;\;\;j\geq2.
\ee
One can check that the last sequence is obtained from the Leja 
sequence $F$ which is constructed recursively by $F_1=\{1\}$
and $F_{2^{n+1}}=F_{2^n}\wedge e^{\frac {i\pi}{2^n}} \o{F_{2^n}}$. Both
$\Re$-Leja sequences $R$ satisfies $2r_0^2-1=1$, $2r_2^2-1=-1$ and more 
generally $2r_{2j}^2-1 = 2r_{2j-1}^2-1 = r_j$
for any $j\geq2$, thanks to the trigonometric identity 
$2\cos^2(\theta/2)-1=\cos(\theta)$. This shows that in both cases 
$R$ satisfies 
the property 
\begin{equation}
R^2=R
\quad \mbox{where}\quad
R^2 := (2r_{2j}^2-1)_{j\geq0},
\label{defR2}
\end{equation}

In general, given a Leja sequence $E$ in $\cU$ initiated at $1$ 
and $R$ the associated $\Re$-Leja sequence, we have that
$R^2$ is an $\Re$-Leja sequence and it is associated 
with $E^2$ which, in view of Theorem 
\ref{TheoImplications}, is also a Leja sequence initiated 
at $1$. This result is given in \cite[Lemma 3.4]{Ch2}
and it has many useful implications that we have exploited 
in order to prove that $\D_{k}(R)$ grows at 
worse quadratically. 

For all Leja sequences $E$ on $\cU$ initiated at $1$, the section 
$E_{2^{n+1}}$ is equal in the set sense with the set of $2^{n+1}$-roots 
of unity, therefore for all $\Re$-Leja sequences $R$, the section 
$R_{2^n+1}$ is equal to the set of Gauss-Lobatto 
abscissas of order $2^n$, i.e. 
\be
\label{GaussLobatto}
R_{2^n+1} = \Big\{\cos(j\pi/2^n): j=0,\dots,2^n\Big\},
\ee 
in the set sense. This set of abscissas is optimal as far as 
Lebesgue constant is concerned, in the sense 
$\L_{R_{2^n+1}}\simeq\frac {2\log (2^n+1)}\pi$.
More precisely, we have the bound
\begin{equation}
\label{lebesgueGaussLobatto}
\L_{R_{2^n+1}}\leq 1+\frac 2 \pi \log(2^n),
\end{equation}
see \cite[Formulas 5 and 13]{DzI}. This suggests 
that the sequence $R$ might have a moderate growth of the 
Lebesgue constant of its section $R_k$. 

In the paper \cite{CaP2}, it has been proved that 
$\L_{R_k} = {\cal O}(k^3\log (k))$. We have improved this 
bound in \cite{Ch1,Ch2} and showed that 
$\L_{R_k} \leq 8\sqrt 2~k^2$ for any $k\geq2$. 
Here we again exploit our approach of \cite{Ch2} 
which, using simple calculatory arguments, 
relate the analysis of the Lebesgue function 
associated with $R_k$ to that of the Lebesgue function 
associated with the smaller Leja section that yields 
$R_k$ by projection. This approach allows us to circumvent 
cumbersome real trigonometric functions which arise in the 
study $\lambda_{R_k}$, see \cite{CaP2,Ch1}, and to 
take full benefit from the machinery developed for Leja sequence on $\cU$.

\begin{remark}
Without loss of generality, we assume for the remainder 
of this section that $E$ is the simple Leja sequence in 
\iref{SimpleLejaSequence} and $R$ the associated $\Re$-Leja
sequence. All our arguments hold in the more general case, 
the assumption is essentially for notational clearness.
It allows us, in view of \iref{memoryShort}, to use $E$ instead for 
$E^2$ and more generally instead of $E^{2^p}$ which is defined by 
$E^{2^p}:=((e_{2^p j})^{2^p})_{j\geq0}$. 
\end{remark}

The bound \iref{lebesgueGaussLobatto} is 
sharp and we are only interested in bounding 
$\L_{R_k}$ when $k-1$ is not a power of $2$. 
For the remainder of this section, 
we use the notation
\be
\label{valuesNKKprime}
\begin{array}{l}
n\geq 0,
\quad\quad 2^n < k-1 < 2^{n+1},
\quad\quad 0< l := k-(2^n+1) < 2^n\\
\quad K:=2^{n+1}+l,
\quad\quad\quad G_k:= E_K,
\quad\quad\quad F_K :=E_{2^{n+1},K}.
\end{array}
\ee  
We should note that  in \cite{Ch2} we have used 
$k'$ and $F_k$ to denote $l$ and $F_K$.
In view of \iref{formulaJ},
we have $K=J(k)$, so that $E_K$ is the smallest 
section that yields $R_k$ by projection into $[-1,1]$. 
We denote by $L_0,L_1,L_2,\cdots,L_{K-1}$ the 
Lagrange polynomials associated with $E_K$. 
The inspection of the the proof of 
\cite[Lemma 6]{Ch2} shows that 
for $z\in \partial \cU$ and $x=\Re(z)$,
\be
\label{boundLRk}
\lambda_{R_k}(x) \leq
\gamma_K(z)
+
\gamma_K(\o z),\quad\quad 
\gamma_K(z):=
|w_{F_K}(\o{z})|
\sum_{j=0}^{K-1} 
\frac {|L_j(z)|}
{|w_{F_K}(\o{e_j})|}.
\ee
In the proof of \cite[Lemma 6]{Ch2}, we have bounded 
the functions $|w_{F_k}|/{|w_{F_k}(\o{e_j})|}$ in the 
previous sum by $2^{n+\frac 12 -p(l)}$. This implied the result of 
\cite[Theorem 5]{Ch2}, namely
$\L_{R_k}\leq 2^{n+\frac 32 -p(l)} \L_{E_K}$.
In view of the new bound \iref{bestBound} 
and the facts that
$p(K)=p(l)$, $\sigma_1(K)=1+\sigma_1(l)$
and $K=2^n+k-1 \leq 3 \times2^n$, the 
previous bound implies 

\be
\label{boundGood}
\L_{R_k} 
\leq
12\sqrt {3}
~
2^{\frac {3n-3p(l)+\sigma_1(l)} 2}
~\L_{2^{p(l)}},\quad k\geq1,
\ee
where $\L_{2^p}$ is bounded as in \iref{boundL2p}.
We propose to improve slightly the previous 
inequality by applying rather Cauchy Schwartz 
inequality when bounding the function $\gamma_K$. 

\begin{theorem}
Let $R$ be an $\Re$-Leja sequence and $n,~k$ and $l$ 
as in \iref{valuesNKKprime}. We have
\be
\label{boundBestRLeja}
\L_{R_k} \leq 6\sqrt{5}~
2^{n+\sigma_1(l)-p(l)} \L_{2^{p(l)}}, 
\ee
where $\L_{2^p}$ is bounded as in \iref{boundL2p}.
\label{theoLebesgueConstant}
\end{theorem}
{\bf Proof:}
In order to lighten the notation, we use the shorthand 
$p$ in order to denote $p(l)$. We introduce $l'$ and 
$K'$ and $F_{K'}$ defined by 
$$
l' := l/2^p, 
\quad \quad
K' :=K/2^p=2^{n-p+1}+l',
\quad\quad
F_{K'} := E_{2^{n-p+1},K'}.
$$
The sequence $E$ satisfies $E^2=E$ and one 
can check that 
$w_{F_K}(z) = w_{F_{K'}}(z^{2^p})$.
Also by $ e_{2j}^2 = e_{2j+1}^2  =e_j$, one 
has $(e_{2^p j+q})^{2^p}=e_j$ for any 
$q=0,\dots,2^p-1$. Moreover, if 
$M_1,\dots,M_{K'-1}$ are the Lagrange polynomials 
associated with $E_{K'}$, then
$$
\sum_{q=0}^{2^p-1}|L_{2^pj+q}(z)| 
\leq \L_{2^p} 
M_j(z^{2^p}),\quad\quad j=0,\dots,K'-1,
$$
see the proof of \cite[Theorem 2.8]{Ch1}.
Therefore by pairing the indices in the sum 
giving $\gamma_{K}$ by $2^p j+q$ for 
$j=0,\dots,K'-1$ and $q=0,\dots,2^p-1$, we infer 
$$
\gamma_K(z)
\leq 
\(
|w_{F_{K'}}(\o{\xi})|
\sum_{j=0}^{K'-1} 
\frac {|M_j(\xi)|}
{|w_{F_{K'}}(\o{e_j})|}
\)
\L_{2^p}
=
\L_{2^p}
\gamma_{K'}(\xi),\quad\mbox{with}\quad \xi=z^{2^p}.
$$
In view of \iref{boundLRk}, this implies that 
$\L_{R_k} \leq 2  \L_{2^p} 
\sup_{\xi\in\cU}\gamma_{K'}(\xi)$.
Applying Cauchy Schwatrz inequality to $\gamma_{K'}$
and using that $F_{K'}$ is an $l'$-Leja sequence, we have 
for any $\xi\in\cU$
$$
\gamma_{K'}(\xi) 
\leq 
2^{\sigma_1(l')}
\(\sum_{j=0}^{K'-1}
\frac 1{|w_{F_{K'}}(\o{e_j})|^2}\)^{1/2}
\(\sum_{j=0}^{K'-1} |M_j(\xi)|^2\)^{1/2}
= 
2^{\sigma_1(l')+n-p} 
\sqrt{\gamma_{n-p+1,l'}}~\lambda_{E_{K'},2}(\xi),
$$
where $\gamma_{n-p+1,l'}$ is defined as in \iref{defBetaln} 
with $m=n-p+1$ and $\lambda_{E_{K'},2}$ is the quadratic 
Lebesgue function associated with $E_{K'}$. In view of the 
bounds we have 
for these quantities and in view of $\sigma_1(K')=1+\sigma_1(l')$ 
and $\sigma_1(l')=\sigma_1(l)$, we get
$$
\gamma_{K'}(\xi)
\leq
2^{\sigma_1(l)+n-p} 
\sqrt{\frac {5}{2^{\sigma_1(l')+1}}}
~3\sqrt {2^{1+\sigma_1(l')}-1}
\leq 
3\sqrt5 
~
2^{\sigma_1(l)+ n- p}. 
$$
The proof is then complete.\hfill $\blacksquare$
\nl

The bound in \iref{boundBestRLeja} improves the bound in 
\iref{boundGood} by $2^{\frac{\sigma_1(l)+p(l) -n}2}$. The 
bound can also yield linear estimates for $\L_{R_k}$, for 
instance when $l$ is such that 
$ 2^{\sigma_1(l)-p(l)} \L_{2^{p(l)}} \leq 1$, which is the case 
if for example $p(l) \geq 2 \sigma_1(l)$.  
However, if $0<l<2^n$ is the integer with 
the most number of ones in the binary expansion, 
i.e. $\sigma_1(l)=n$ or $l=2^n-1$ and 
$k=2^{n+1}$, we merely 
get the quadratic bound 
\be
\label{pessimisticQuadratic}
\L_{R_k} \leq 6\sqrt 5~ 2^{2n}= \frac {3\sqrt 5}2 k^2.
\ee
In \cite{CaP2}, section 3.4, it is shown that for the values 
$k=2^{n}$, in other words $R_k$ is the set of Gauss-Lobatto 
abscissas \iref{GaussLobatto} missing one abscissa, one has 
$\L_{R_k}\geq \lambda_{R_k}(r_k)= k-1$. As a consequence,
the growth of $\L_{R_k}$ for $k\geq1$ can not be
slower than $k$. However, for this case, we can prove
$\L_{R_k}\leq 3 k$, see \iref{finalRemark}, showing that 
\iref{pessimisticQuadratic} is rather pessimistic.

The estimate in \iref{boundBestRLeja} is logarithmic for 
many values of the integer $k$. For instance, if
$k=(2^n+1)+2^{n-p} k'$ 
for some $p=1,\dots,n$ and some $k'=0,\dots,2^p-1$, then we have 
$l=2^{n-p} k'$, so that $n-p \leq p(l)\leq n$ and 
$\sigma_1(l)=\sigma_1(k')\leq p$ implying that
\be
\L_{R_k} 
\leq 
6 \sqrt 5 ~2^{2p}~ \L_{2^{p(l)}}
\leq
6 \sqrt 5 ~2^{2p}~ \L_{2^n}
\leq 
6 \sqrt 5 ~2^{2p} 
\frac 2\pi
\(\log(2^n)+ 9/4\).
\ee
For a small value of $p$, the previous estimate is as good 
as the optimal logarithmic estimate $\frac {2\log(k)}\pi$ for 
large values of $n$. Given then $p$ fixed, one has $2^p$ 
intermediate values between $2^n+1$ and $2^{n+1}+1$, 
which are the numbers $k=(2^n+1)+2^{n-p} k'$ for 
$k'=0,\dots,2^p-1$, for which the Lebesgue constant is 
logarithmic. This observation can be used in order to 
modify the doubling rule with Clemshaw-Curtis abscissas 
in the framework of sparse grids, see \cite{GWG}. 

\section{Growth of the norms of the difference operators}
\label{sectionDiff}

In this section, we discuss the growth of the norms of the 
difference operators $\Delta_0 = I_{Z_1}$ and 
$\Delta_{k} = I_{Z_{k+1}} - I_{Z_k}$ for $k\geq1$, associated 
with interpolation on Leja or $\Re$-Leja sequences. We are 
interested in estimating their norms $\D_k$ defined in 
\iref{normDeltak}.
Elementary arguments, see \cite{Ch2}, show that
\be
\label{normDiff}
\mathbb D_k (Z)= 
\Big(1+\lambda_{Z_k} (z_k) \Big)
\sup_{z\in X} \frac {|w_{Z_k}(z)|}{|w_{Z_k}(z_k)|},\quad k\geq1.
\ee
In particular if $Z$ is a Leja sequence on the compact $X$, then 
\be
\D_k(Z)=1+\lambda_{Z_k} (z_k).
\ee 
In \cite{Ch1}, we have established that 
$\lambda_{E_k} (e_k)\leq k$ if $E$ is a Leja 
sequence on $\cU$ initiated at $\partial \cU$,
which implies $\D_k(E)\leq1+k$. Here, we improve 
slightly this bound. As for the improvement of 
\iref{boundLkLatest} into \iref{bestBound}, we have
\begin{theorem}
Let $E$ be a Leja sequence on the unit 
disk initiated at $e_0\in \partial \cU$, One has 
$\D_0(E)=1$ and 
\be
\D_k(E) \leq 
1+\sqrt {\frac {k}{2^{p(k)}}(2^{\sigma_1(k)}-1)}
~~\L_{2^{p(k)}}
\ee
\end{theorem}

For $\Re$-Leja sequences $R$ on $[-1,1]$, we have shown 
in \cite{Ch2} using a recursion argument based on the fact 
that $R^2$ defined as in \iref{defR2} is also an $\Re$-Leja 
sequence, that   
\be 
\D_k(R) \leq (1+k)^2,\quad\quad k\geq0.
\ee
In view of the new bounds obtained in this paper for Lebesgue 
constant of $\Re$-Leja sections, the previous bound is not 
sharp. Indeed, we have 
$\D_k\leq\L_k+\L_{k-1} \leq 12\sqrt{5}~ k^{3/2}$,
for $k$ such that $l=k-(2^n+1)\leq 2^{n/2}$. We give 
here a sharper bound for $\D_k(R)$. 
We recall that up 
to a rearrangement in the formula \iref{normDiff}, see \cite{Ch2} for 
justification, we may write the quantities $\D_k(R)$ in a more 
convenient form for $\Re$-Leja sequences. We introduce the 
polynomial $W_{R_k} := 2^k w_{R_k}$, we have
\begin{equation}
\label{deltabeta}
\mathbb D_k(R) = 
2
\beta_k(R)
\sup_{x\in [-1,1]}  |W_{R_k}(x)|,
\;\;
\;\;
\beta_k(R):=
\frac {1+\lambda_{R_k}(r_k)}{2|W_{R_k}(r_k)|},
\end{equation}
We have already proved in \cite[Lemma 7]{Ch2} that
\be
\label{beta2nk}
\beta_{2^n}(R)=1/4
\quad\mbox{and}\quad 
\beta_k(R) \leq 2^{\sigma_0(k)-p(k)-1},
\quad\mbox{for}\quad k\neq 2^n.
\ee
Here we provide a 
sharper bound for $\D_k(R)$ by slightly improving 
the estimate $4^{\sigma_1(k)+p(k)-1}$ that we have 
established in \cite{Ch2} for
$\sup_{x\in [-1,1]}  |W_{R_k}(x)|$.
\begin{lemma}
Let $R$ be an $\Re$-Leja sequence in $[-1,1]$, $n\geq1$, 
$2^n+1\leq k<2^{n+1}+1$ and $l=k-(2^n+1)$.
One has $\sup_{x\in [-1,1]}  |W_{R_{k}}(x)|\leq 2^{n+3}$ if 
$k=2^{n+1}$, else
\be
\sup_{x\in [-1,1]}  |W_{R_k}(x)|
\leq 
2^{2\sigma_1(k)+p(k)-1}.
\ee
\end{lemma}
{\bf Proof:}
We use the notation $K$, $G_k$ and $F_K$ as in 
\iref{valuesNKKprime} and introduce $G_{k+1}:=E_{K+1}$ 
and $F_{K+1}:=E_{2^{n+1},K+1}$.
In view of \cite[Lemma 5]{Ch2}, one has for 
$z\in\partial \cU$ and $x=\Re(z)$
$$
|W_{R_k}(x)|
=
|z^2-1||w_{G_k}(z)| |w_{F_K}(\o z)|
=
|z-\o z||w_{G_k}(z)| |w_{F_K}(\o z)|.
$$
Also since $|z-\overline{z}|\leq |z-e_K|+|\o z-e_K|$, then
$$
|W_{R_k}(x)| \leq
|w_{G_{k+1}}(z)||w_{F_K}(\o z)|
+|w_{G_k}(z)||w_{F_{K+1}}(\o z)|.
$$
In the two previous inequalities, one has $F_K=\emptyset$ and 
$w_{F_K}\equiv1$ in the case $k=2^n+1$.
 We have that $G_k$, 
$G_{k+1}$, $F_K$ and $F_{K+1}$ are all Leja sections with length 
$K,~K+1,~l$ and $l+1$ respectively. Therefore, by the 
second property in Theorem \ref{TheoImplications}
$$
|W_{R_k}(x)| 
\leq 
\min\(2^{1+\sigma_1(K)+\sigma_1(l)}
,
2^{\sigma_1(K+1)+\sigma_1(l)}
+
2^{\sigma_1(K)+\sigma_1(l+1)}
\)
=2^{2+\sigma_1(l)}
\min(2^{\sigma_1(l)},2^{\sigma_1(l+1)}),
$$
where we have used 
$\sigma_1(K)=1+\sigma_1(l)$ and 
$\sigma_1(K+1)=1+\sigma_1(l+1)$
since $K=2^{n+1}+l$ and $l<2^n$.
If $k=2^{n+1}$, i.e. $l=2^n-1$, then
$|W_{R_k}(x)|\leq 2^{3+n}$. Else
by $k=2^n+(l+1)$ and $0\leq l<2^n-1$,
$$
\sigma_1(k) -1= \sigma_1(l+1)\quad\mbox{and}\quad
\sigma_1(k)-2+p(k)=\sigma_1(k-1)-1 = \sigma_1(l).
$$
Therefore 
$$
|W_{R_k}(x)|
\leq 2^{2\sigma_1(k)+p(k)-1}
\min(2^{-1+p(k)},1),
$$
which completes the proof.
\hfill $\blacksquare$
\nl

By injecting the estimate of the previous lemma 
and the estimate of \iref{beta2nk} in formula 
\iref{deltabeta} and by using the identity 
$\sigma_0(k)+\sigma_1(k)=n+1$ for 
$2^n\leq k<2^{n+1}$, we are able to conclude 
the following result.
\begin{cor}
Let $R$ be an $\Re$-Leja sequence in $[-1,1]$.
The norms of the difference operators 
satisfy, $\mathbb D_0(R) =1$ and 
for $2^n \leq k<2^{n+1}$
\begin{equation}
\mathbb D_k(R) \leq 2^{\sigma_1(k)} 2^{n}
\end{equation}
\end{cor}

\begin{figure}[!h]
\includegraphics[scale=0.6]{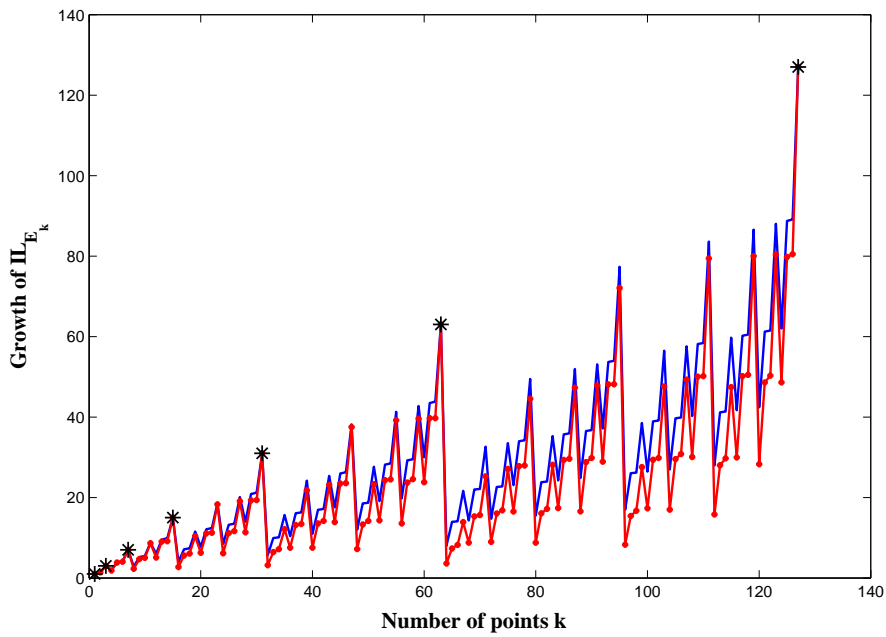}\\
\includegraphics[scale=0.6]{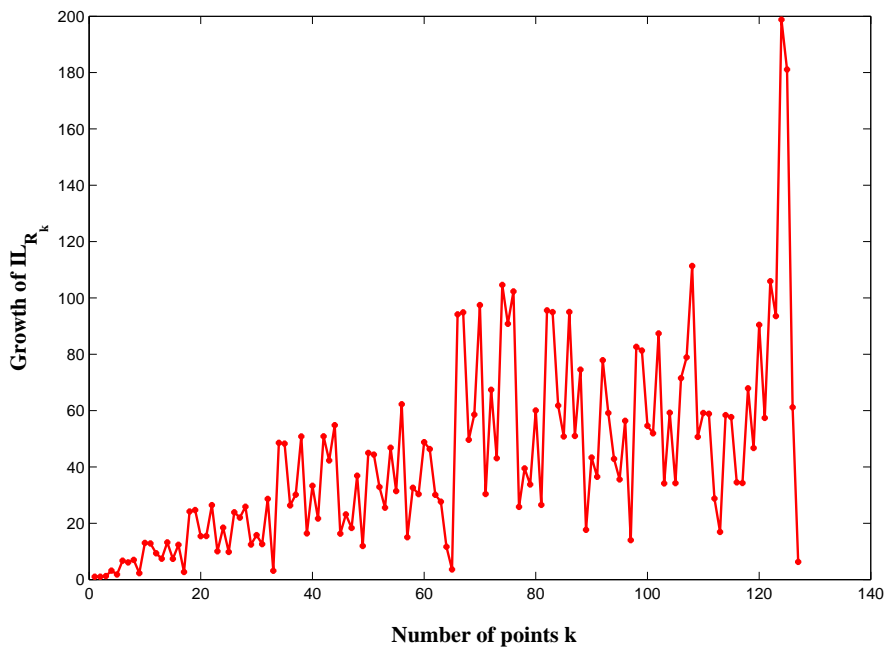}
\caption{Exact Lebesgue constants
associated to the $k$-sections of the Leja sequence 
$E$ and the assciated $\Re$-Leja sequence $R$ for 
$k=1,3,\dots,129$.}
\label{fig}
\end{figure}

The previous estimates can be used in order 
to provide estimates for $\L_{R_k}$ that can be sharper 
than \iref{boundBestRLeja}. We have $\Delta_k=I_k - I_{k-1}$, 
therefore
\be
|\L_{R_{k+1}}-\L_{R_k}|\leq \D_k(R),\quad k\geq1.
\ee
In particular, the estimate in the previous corollary 
combined with the sharp bound
\iref{lebesgueGaussLobatto} implies that for the value 
$k=2^n$, we get 
\be
\L_{R_k} 
\leq 1+\frac 2 \pi \log(2^n)+ 2^{n+1}
\leq 
3k
\ee
This shows that in the case $k=2^n$ which corresponds to 
$R_k$ being the set of Gauss-Lobatto abscissas \iref{GaussLobatto} missing one abscissas and 
for which $\L_{R_k}\geq k$, the previous bound 
is satisfactory. This also confirm that the estimates 
\iref{boundBestRLeja} is indeed pessimistic in this 
case, see the inequality \iref{pessimisticQuadratic}.
This added to the observed growth of $\L_{R_k}$ 
for values $k\leq 128$, Figure \ref{fig}, suggests 
that the bound 
\be
\label{finalRemark}
\L_{R_k} \leq 3k,\quad\quad k\geq1,
\ee 
might be valid for any $\Re$-Leja sequence $R$. We 
conjecture its validity.

In Figure \ref{fig}, we also represent for the values
$k\leq 128$, the growth of the Lebesgue constant 
$\L_{E_k}$ (in blue) and the estimate $\sqrt {k(2^{\sigma_1(k)}-1)}$
(in red) which multiplied by $3$ bounds $\L_{E_k}$, see
\iref{boundLkLatest}. We observe that the regular patterns in the graph 
of $k \mapsto \L_{E_k}$, which reveals the particular role of divisibility 
by powers of 2 in $k$, is caught by the estimate. The worst values of 
$\L_{E_k}$ appear for the values $k = 2^n-1$ for which it was proved in 
\cite{CaP1} that $\L_{E_k} = k$ and which is also equal to 
$\sqrt {k(2^{\sigma_1(k)}-1)}$ since $\sigma_1(k)=n$.

\newpage


\ifx\undefined\bysame
\newcommand{\bysame}{\leavevmode\hbox to3em{\hrulefill}\,}
\fi

\end{document}